\documentclass[conference]{IEEEtran}
\IEEEoverridecommandlockouts
\usepackage{cite}
\usepackage{amsmath,amssymb,amsfonts}
\usepackage{algorithmic}
\usepackage[caption=false,font=footnotesize]{subfig}
\usepackage{graphicx}
\usepackage{textcomp}
\usepackage{xcolor}
\usepackage{multicol}
\usepackage[ruled,vlined]{algorithm2e}
\usepackage{geometry}
\def\BibTeX{{\rm B\kern-.05em{\sc i\kern-.025em b}\kern-.08em
    T\kern-.1667em\lower.7ex\hbox{E}\kern-.125emX}}
\begin{document}

\title{Sphereabout: A Spherical Intersection Design for Networked Urban Air Mobility\\
\vspace{1cm}
\small
IEEE Intelligent Transportation Systems Conference, 2025, Gold Coast, Australia
\thanks{This research is funded by the Canada Research Chair in Disruptive Transportation Technologies and Services (CRC-2022-00480).}
}

\author{\IEEEauthorblockN{Farzan Moosavi}
\IEEEauthorblockA{\textit{Laboratory of Innovations In Transportation} \\
\textit{Torono Metropolitan University}\\
Toronto, Canada \\
0009-0008-8490-3851}
\and
\IEEEauthorblockN{ Bilal Farooq}
\IEEEauthorblockA{\textit{Laboratory of Innovations In Transportation} \\
\textit{Torono Metropolitan University}\\
Toronto, Canada \\
0000-0003-1980-5645}
}

\maketitle

\begin{abstract}

Urban aerial mobility is rapidly expanding, specifically on-demand Unmanned Aerial Vehicle (UAV) delivery services in urban environments. This necessitates management of the low-altitude airspace network to ensure smooth and safe traffic throughput. This study introduces a novel three-dimensional aerial network design inspired by a terrestrial transportation graph network for UAV networked mobility. Utilizing two-way tube corridors for node-to-node delivery and a spherical roundabout model, as ``Sphereabout'', to optimize the traffic flow through the spherical intersection. Through this architecture, three-dimensional conflict-free air mobility management can be achieved via numerical experiments and utilizing the geometrical features of this intersection.

\end{abstract}


\section{Introduction}

A well-defined aerial network structure is required for safety, efficiency, and capacity management of Unmanned Aerial Vehicles (UAVs) in high-density urban areas. Specifically, urban air mobility in unconstrained airspace requires an operational and realistic delivery network, traffic alignment and segmentation to mitigate collision probability. 
For example, point-to-point delivery through a straight line path might lead to obstacle collision, such as buildings. As a result, virtual air corridors for point-to-point safe passage of multiple UAVs are represented through a Voronoi over-building diagram by a well-clear separation from buildings \cite{b3}. Moreover, a road-based network that automatically avoids conflict and a collision-free multi-layered discretization of airspace are introduced \cite{b4}. The latter has more freedom for UAVs, freely choosing their position, altitude, heading, and speed, which increases airspace capacity and reduces flying costs. However, conflict resolution in such a low-altitude urban environment is essential. For instance, a risk-based airspace is considered to demonstrate adaptive decision-making strategies for different flight conflicts to optimize the conflict resolution of unmanned aircraft systems \cite{b5}. Although their proposed algorithm is effective under different traffic density scenarios, it cannot perform well when the traffic density increases significantly. Thus, the scalability problem of the optimization algorithm still has room for improvement. Similarly, \cite{b6} presented Unmanned Aircraft System Traffic Management (UTM) with conflict detection and resolution (CD\&R), illustrated a linearized version of a non-linear complete optimization model to solve a higher scale problem, having a slight compromise on solution optimality.

\begin{figure}[htbp]
\centerline{\includegraphics[width=0.4\textwidth]{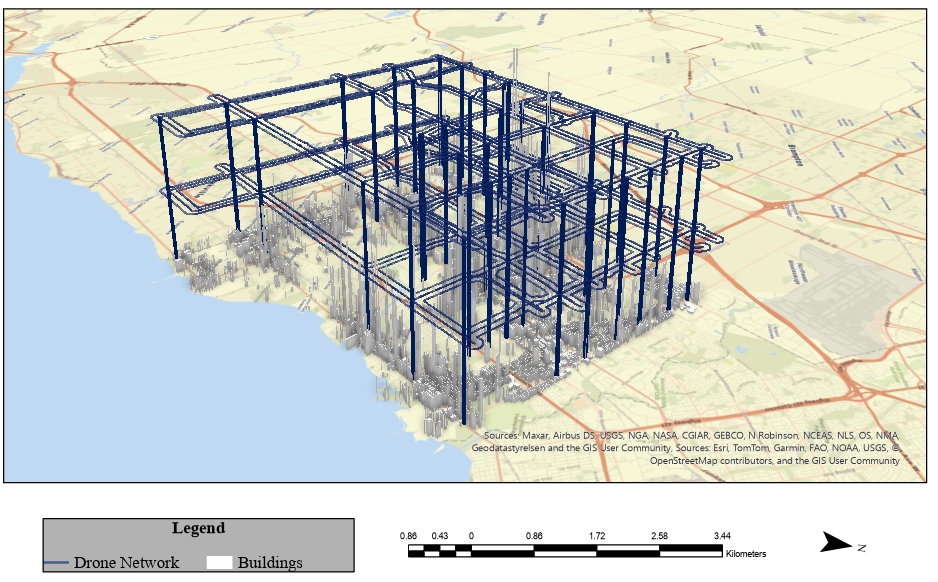}}
\caption{UAV network with elevated virtual links on top of existing roads.}
\label{f1}
\end{figure}

On the other hand, a highly structured network like the elevated ground transportation network would lead to a safer and obstacle-free path planning \cite{b2}. In fact, most works agree with the hybrid concept that preventing conflicts is better than resolving conflicts \cite{b7}. Also, to scale and generalize the airspace design for managing urban mobility, a flexible airspace must be adapted to city characteristics and regulations. In this regard, the over-road network with multiple layers, as shown in Fig. \ref{f1}, is an effective solution if potential conflict is resolved at intersections. Considering the three-dimensional aerial network, a new type of intersection must be designed to regulate the horizontal, lateral, and vertical transitions so UAVs can change their path into different road links or altitudes. 

This work introduces Sphereabout, a spherical roundabout intersection to accommodate UAV movement in networked urban air mobility (Fig. \ref{f2}). Six directions are considered two-way tube corridors with a buffer zone, leading to six entry tubes and six exit tubes in the opposite direction. Nevertheless, several concepts of the three-dimensional intersection model are proposed \cite{b10, b11}. For instance, introducing a space-time capsule as a protected area, a UAV can join or leave freely into the tube corridor and the safety lane.
Additionally, a combination of turning rules in a roundabout connected to vertical corridors with ascending and descending zones to climb. However, the spherical design can improve traffic flow using the alternative path on the sphere. Besides, to the best of the author's knowledge, no research has studied the traffic flow assignment for three-dimensional spherical intersections. 
In summary, our contribution is as follows:

\begin{itemize}
    \item Novel spherical intersection design and operation for air mobility network inspired by ground roundabout with arc-based transportation modes.
    \item A conflict-free path assignment ensuring smooth and safe traffic flow.
    \item An optimization model for maximum traffic throughput at the intersection and traffic flow management. 
\end{itemize}

\section{Background} \label{back}

In the highly dense urban environment, the most critical factor is to ensure safety while maintaining optimal traffic flow at links and intersections. Existing research has proposed obstacle-free vehicle routing in urban areas and UAV regulations, such as compliance with air regulations, including no-fly zones, height restrictions, and safety protocols \cite{b8}; however, less structured airspace allows for higher traffic densities and lacks scalability to the system, and too much structure reduces performance as flight paths become constrained and capacity efficiency decreases. 
Accordingly, the trade-off has been made by \cite{b12} for various approaches to airspace design. Based on this ranking, the notion of a structured network for aerial delivery is created by extending the current 2D road network to the third dimension to capture the inheritance of collision avoidance of the city structure and bringing traffic control at the aerial intersection. To this end, it is critical to develop an efficient traffic management system at these aerial intersections to accommodate through and side turn traffic, as well as vertical turn feasibly and free of collisions.


To this end, \cite{b9} investigated two novel airspace design concepts, the two-way and one-way concepts, for the constrained urban environment. Based on traffic alignment and segmentation principles, both airspace concepts employed heading-altitude rules to separate cruising traffic for their heading directions vertically. Eventually, they demonstrate that in a constrained urban environment, having vertically segmented altitude layers to accommodate traffic with similar directions and some horizontal constraints imposed on traffic flow is beneficial for safety. Despite their systematic rule-based traffic segmentation, conflict resolution at the intersection was not studied, and it incurs an excessive transportation cost since the choice of changing lanes or headings will need to be taken at a specific altitude and not at any intersection. Moreover, \cite{b15} presented air mobility intersection traffic control using a Multi-Agent Reinforcement Learning (MARL) framework to manage the two-dimensional traffic flow by accelerating and decelerating UAVs to avoid conflict.
Furthermore, \cite{b16, b17} proposed multi-lane intersection planning for drone corridors. The first introduced a durational speed scheduling rather than lane change and hovering to predict potential conflicts and compute new velocities based on relative distance and heading angles. Their simulation-based algorithm detects the possible conflict and applies a velocity change during the specific time window to pass the intersection at different times without violating the minimum separation. The second work presented a corridor changing planning based on the lane connectivity to prevent conflicts while minimizing deviation and delay. Various maneuvers simulations, such as three-dimensional path planning across lanes for detecting conflict, are considered for path modification, leading to maintaining safe separation. Nonetheless, neither work considered the vertical lane change in the sense of layer transition to the higher and lower altitudes. Similarly, research \cite{b18} demonstrated a three-dimensional reservation-based scheduling system for intersection traffic management by optimizing the sequence in which UAVs enter the intersection using a Genetic algorithm. However, the airspace intersection aims to utilize the capacity from all the traffic flow directions, including vertical flows like transition from the third dimension to the planar network. Therefore, in this work, a conflict-free traffic flow optimization is considered to obtain an efficient path planning.

\section{Problem Description}\label{prob}


\subsection{Spherical Intersection Design}\label{AA}

First, the overall diameter of the sphere is considered no less than that of a typical street in the city of Toronto with a maximum speed of 40 $km/h$, which is 26 meters. This includes the two-way street, parking, bicycle lane, and curbside width. As a result, the sphere radius, $R$, is 13 meters. Initially, the local coordinate is laid out as if the x-direction is along the East side and the y-direction is along the North. Consequently, the z-direction would align with the radial direction, pointing to the zenith.  

One significant aspect of spherical design is to utilize the symmetry in geometry. Therefore, as there are eight tubes perpendicular to the z-direction, the sphere is sliced into eight equal segments, and each tube is placed at an equal angular distance, facing the corresponding direction. In addition, each tube diameter is set to ensure enough safety distance among UAVs to avoid the interference introduced by downwash effects \cite{b20}. Also, it depends on the velocity of the UAV; for a forward flight, the x-axis clearance needs to be increased. This clearance distance affects the horizontal, lateral, and vertical stream differently. For example, a quadcopter, with a rotor diameter, $D$, would need a lateral and vertical separation gap as in equations \eqref{eq1} and \eqref{eq2}, respectively \cite{b21}. 

\begin{equation}
 x \geq 4D \label{eq1}, \;\; V = 5\;m/s
\end{equation}  
\begin{equation}
 z \geq 1.5D \label{eq2}, \;\; V = 5\;m/s
\end{equation}  

\begin{figure*}
     \centering
     \subfloat[3D view.]{\includegraphics[width=0.22\textwidth]{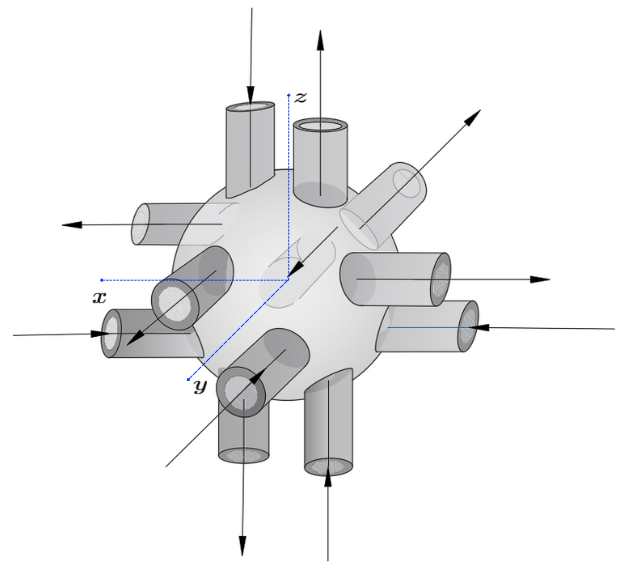}\label{f2}}
     \hfill
      \subfloat[x-y view.]{\includegraphics[width=0.2\textwidth]{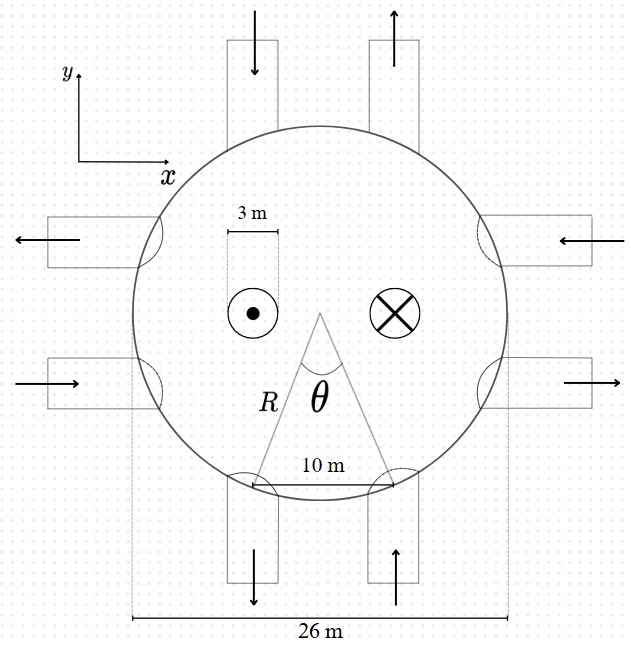}\label{f3}}
     \hfill
      \subfloat[x-z view.]{\includegraphics[width=0.2\textwidth]{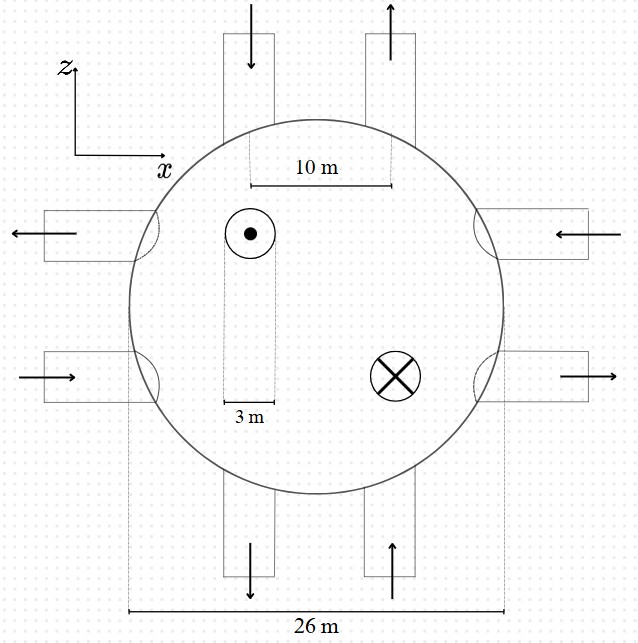}\label{f4}}
     \hfill
      \subfloat[z-y view.]{\includegraphics[width=0.2\textwidth]{xz.jpg}\label{f5}}      
      \caption{Sphereabout design configuration and the corridor flow directions.}
\end{figure*}

Where $V$ is the velocity of the UAV, therefore, for the DJI FlyCart30, the rotor diameter is 1.375 meters, the inner tube keeping away from downwash is considered 2 meters with an additional 1 meter as a buffer radial zone. The spherical sketch from the principal axis view, front, lateral, and top views are depicted in Fig. \ref{f3}, \ref{f4}, and \ref{f5}, respectively. Each tube's minimum perpendicular separation, $d$, can be derived by equation \eqref{eq3}, satisfying the minimum clearance distance, ensuring the safety of the operation.


\begin{equation}
  d = Rsin(\frac{\theta}{2}) \approx 10 \;m\label{eq3}
\end{equation}  

\subsection{Path Assignment}

\begin{figure}
     \centering
     {\includegraphics[width=0.5\textwidth]{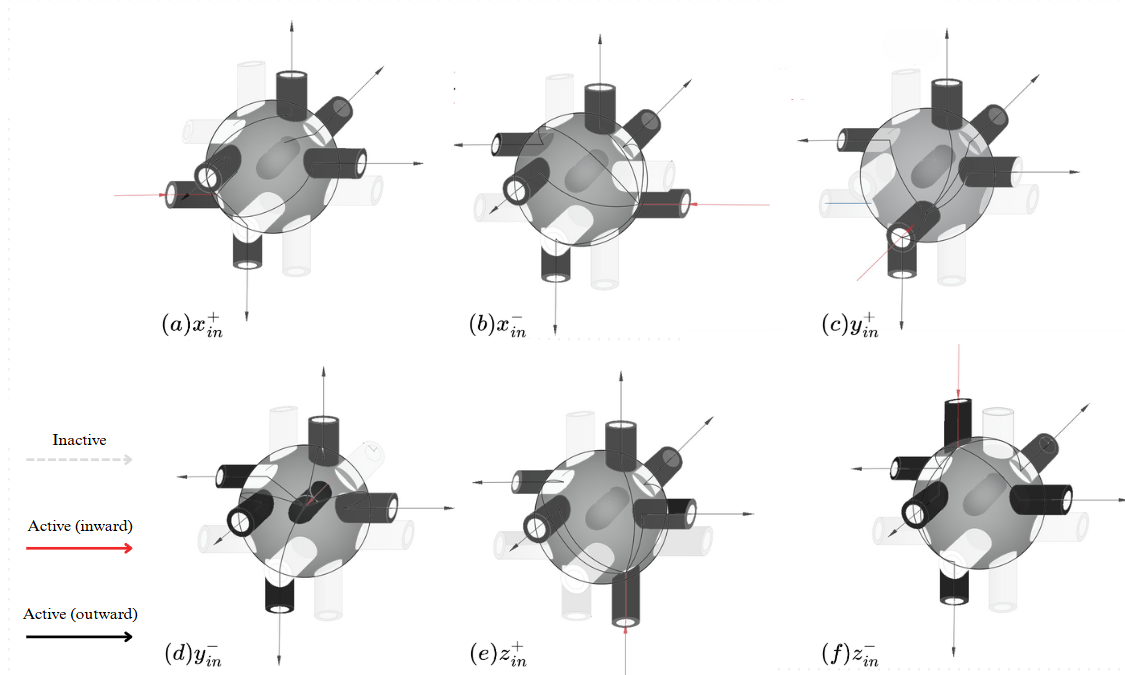}}
       \caption{Entry flows and feasible exit flows.}
\label{f6}
\end{figure}

The advantage of using the spherical design is the alternative of moving along the sphere arc and the straight point-to-point flight. Furthermore, the great circle that passes through two desired nodes can be traversed from two sides: the closest arc between them and the longer arc. This design uses a bipartite graph of six entry nodes for inward flow, $ E = \left\{ x^{+}_{in}, x^{-}_{in}, y^{+}_{in}, y^{-}_{in},z^{+}_{in}, z^{-}_{in}  \right\} $ and six exit nodes for outward flow, $ X = \left\{ x^{+}_{out}, x^{-}_{out}, y^{+}_{out}, y^{-}_{out},z^{+}_{out}, z^{-}_{out}  \right\} $. An illustration of a feasible entry-exit pair along each inward flow direction is shown in Fig. \ref{f6}. Note that the re-entry to the same direction is not allowed, such as from $x^{+}_{in}$ to $x^{-}_{out}$; thereby, each entry node can exactly connect to five exit nodes. Therefore, a bipartite connection from entry to exit, denoted as a set $A=\{(i, j): i \in E, j \in X, \text { no re-entry }\}$.
The alignment of the tubes along the z-axis enables the through traffic without necessarily using the circular path. It can move along the z-axis for the z-direction traffic flow. This can apply to other nodes, though a circular path might be preferred due to the potential conflict. Therefore, all the possible path types a UAV can be assigned from entry to exit node are $P_{ij}  \;\; \forall(i,j) \in A = \{1:\text { direct line, } 2:\text { short arc, } 3:\text { long arc }\}$. In fact, for each feasible pair $(i,j)$, there are three geometric alternatives. 


  


\subsection{Collision Detection}

The main criterion to avoid collision is to ensure the safe clearance distance of the UAV buffer zone is maintained in the Sphereabout, which was set as $d_{min}$ as a dynamic geofence which moves attached to the UAV. As a result, considering both UAV approaching in the direction in which this distance may be violated, the safety distance is regarded as $d_{min} = 3\; m$. Moreover, as in equation \eqref{eq5}, the conflict parameter, $\delta_{kl}$ is defined to determine if for any two distinct paths $m=(i, j, p)$ and $n=\left(i^{\prime}, j^{\prime}, p^{\prime}\right)$ the collision is detected.  

\begin{equation}
\delta_{mn}= \begin{cases}1, & \min _{t, s}\left\|\mathbf{r}_k(t)-\mathbf{r}_{\ell}(s)\right\| \leq d_{min} \\ 0, & \text { otherwise }\end{cases}\label{eq5}
\end{equation}  

where $\mathbf{r}_k(t)$ and $\mathbf{r}_{\ell}(s)$ are the 3D parameterizations of the two candidate paths. This is a non-linear constraint and can be practically discretized at time intervals to check for different entry times that may end up overlaying the buffer zone. As the paths are pre-determined, it would suffice if the minimum distance could be found to check potential collision and instead use the intrinsic length difference of the considered path to compensate for the potential time delay leading to conflict. This led to substantially simplifying and linearizing the optimization problem. The conflict-aware path planning is depicted in the Algorithm \eqref{al1}.

\begin{algorithm}[h]
\caption{Collision Detection and Path Planning}
\label{al1}
\footnotesize
\DontPrintSemicolon
\KwIn{Entry nodes $E$, exit mapping $X$, path types $P$, safety distance $d_{\min}$, sample resolution $N_{\rm samp}$.}
\KwOut{For each scenario of $N$ UAV: an assignment of $(i,j,k)$ or \texttt{NoSolution}.}

\SetKwFunction{ComputeCoords}{ComputeNodeCoords}
\SetKwFunction{SamplePath}{SamplePath}
\SetKwFunction{DetectConflicts}{DetectConflicts}
\SetKwFunction{AssignPaths}{AssignPaths}

\BlankLine
\tcp{1. Precompute node coordinates}
$\{\mathbf p_n\}\leftarrow$\ComputeCoords{$E_i, X_j$}\;

\BlankLine
\tcp{2. Precompute global path samples}
\ForEach{$(i,j)\in A$, $p\in P$}{
  $P_{i,j,p}\leftarrow$\SamplePath{$\mathbf p_i,\mathbf p_j,p,N_{\rm samp}$}\;
}

\BlankLine
\tcp{3. Global conflict detection}
$\mathcal C_{\rm global}\leftarrow\{\}$\;
\ForEach{distinct paths $(i,j,k)\neq(p,q,\ell)$}{
  \If{$\min_{a\in P_{i,j,k},\,b\in P_{p,q,\ell}}\|\!a-b\!\|<d_{\min}$}{
    Add $\bigl((i,j,k),(p,q,\ell)\bigr)$ to $\mathcal C_{\rm global}$\;
  }
}

\tcp{4. Path planning}
\BlankLine
\SetKwProg{Fn}{Function}{:}{}
\Fn{\AssignPaths{$K,P,\mathcal C$}}{
  \tcp{serve all $N$ UAVs}
  \For{$K=N$ \KwTo $1$}{
    \ForEach{subset $K\subseteq UAV$}{
      \ForEach{path choices $\{p_k\in P\mid k\in K\}$}{
        let \texttt{Served} $\leftarrow \{(k,p_k)\mid k\in K\}$\;
        \If{$ (d,\ell),(e,m)\in \texttt{Served} :\;((d,\ell),(e,m))\notin\mathcal{C}$}{
          \Return \texttt{Served}\;
        }
      }
    }
  }
  \KwRet $\emptyset$\;
}
\end{algorithm}

\subsection{Optimization Framework}

The design objective is to maximize the traffic flow at the spherical intersection by a conflict-free path assignment. For each given $N$ UAVs from the set $k=\{1,2, \ldots, N\}$, the traffic flow is known as entry node $i \in E$ and exit node $j \in X$. One binary decision variable for UAV-path assignment is $f_{i j} \in\{0,1\}=1$ if one UAV travels from entry $i$ to exit $j$ and zero otherwise. The other binary variable is 
$x_{i j}^p \in\{0,1\}$ if the UAV on $(i, j)$ uses path $p \in P_{i,j}$; otherwise it is zero. Traffic flow optimization to maximize the number of UAVs passing through the intersection can be found in equation \eqref{eq6}, followed by problem constraints.

\begin{equation}
\max \sum \sum w_{i j} f_{i j} \;\; \forall (i,j) \in A \label{eq6}
\end{equation}  

It is set $w_{i j}=1$, implying only one UAV is utilized in each traffic flow.  

\begin{equation}
\sum_{p \in P} x_{i j}^p=f_{i j}, \;\; \forall (i,j) \in A \label{eq7}
\end{equation}  

\begin{equation}
x_{i j}^p+x_{i^{\prime} j^{\prime}}^{p^{\prime}} \leq 1,  \quad \delta_{(i, j, p), \left(i^{\prime}, j^{\prime}, p^{\prime}\right)} = 1\label{eq8}
\end{equation}

\begin{equation}
x^p_{i j}, f_{ij} \in\{0,1\}, \quad \forall i,j, p\label{eq9}
\end{equation}

The constraint \eqref{eq7} ensures each UAV selects exactly one path and maintains link flow,  \eqref{eq8} refers to a conflict-free constraint where UAVs cannot select paths precomputed to conflict, and \eqref{eq9} satisfies the binary variable constraint. 

\section{Simulation} \label{simul}
The assumptions for the problem scenarios include, (a) UAV paths are known, and only one UAV per corridor is assigned, (b) UAV speed is constant, (c) All the possible traffic flows can be made as long as they are conflict-free, (d) All UAVs start their entry at the same time. Regarding the last assumption, it may be such that the minimum required separation is violated at some point for two arbitrary paths; however, the delay in the entry time can resolve the potential conflict. Nonetheless, this optimization aims to achieve the intrinsically conflict-free intersection architecture so that the traffic can be managed independently of entry time using path assignment. 

\subsection{Experimental Setting}\label{SCM}

For the simulation, we used classical optimization models using Gurobi solver to evaluate the performance of the optimization model. Different experiments, including various numbers of UAV with random traffic flow, have been set to obtain feasible path assignments at spherical intersections to develop a three-dimensional static traffic management system. The UAV velocity is considered as $ 5 \;\frac{m}{s}$. 

\subsection{Results}

It is desired to identify the collision-free scenarios for any combination of UAV constrained by their traffic flow distributions. Therefore, simulations go through with the least entry tube occupancy, namely, two, and progress to the most critical case where all the entries are occupied. If our conflict resolution fails, the most frequent conflicting pairs need to detect instances where the intersection cannot maximize the traffic flow. This can also be resolved by finding the time difference of arriving at the collision point and adjusting the speed to avoid such a conflict.  
The buffer zone, which determines safety and perturb-free traffic, is fixed. Afterwards, from all feasible sets of scenarios in which the traffic flow is maximum, the solutions of path assignment can be found, consisting of ones that require no conflict resolution and the other that turn into conflict-free traffic after path assignment. The flows that end up in conflict are also detected.  
Table \ref{tab2} demonstrates the Sphereabout intersection throughput comparison result for each UAVs combination based on specific parameters. This Table obtains the number of non-feasible, required conflict resolution, and no-conflict solutions. Moreover, the number of feasible scenarios, the following parameters per scenario, the average optimal flow, and the average path load of the selected type are also calculated.

\begin{table*}[htbp]
\caption{COMPARISON RESULT FOR THE TRAFFIC FLOW OPTIMIZATION}
\footnotesize
\hspace{0.3cm}
\resizebox{\textwidth}{!}{
\centering
\begin{tabular}{|c|cccc|cccc|cc}
\cline{1-9}
\textbf{\begin{tabular}[c]{@{}c@{}}R =  13 m\\ d\_min = 3 m\end{tabular}} & \multicolumn{4}{c|}{\textbf{Number of}}                                                                                                                                 & \multicolumn{4}{c|}{\textbf{Average per Scenario}}                                                                                                                &           &           \\ \cline{1-9}
\textbf{\# of UAV}                                                        & \multicolumn{1}{c|}{\textbf{Scenarios}} & \multicolumn{1}{c|}{\textbf{Collisions}} & \multicolumn{1}{c|}{\textbf{No Conflict Solutions}} & \textbf{Conflict Resolution} & \multicolumn{1}{c|}{\textbf{Average Traffic Flow}} & \multicolumn{1}{c|}{\textbf{Path Load 1}} & \multicolumn{1}{c|}{\textbf{Path Load 2}} & \textbf{Path Load 3} &           &           \\ \cline{1-9}
2                                                                         & \multicolumn{1}{c|}{375}                & \multicolumn{1}{c|}{0}                   & \multicolumn{1}{c|}{314}                            & 61                           & \multicolumn{1}{c|}{2.000}                         & \multicolumn{1}{c|}{1.914}                & \multicolumn{1}{c|}{0.005}                & 0.080                &           &           \\ \cline{1-9}
3                                                                         & \multicolumn{1}{c|}{2500}               & \multicolumn{1}{c|}{0}                   & \multicolumn{1}{c|}{2436}                           & 64                           & \multicolumn{1}{c|}{3.000}                         & \multicolumn{1}{c|}{2.760}                & \multicolumn{1}{c|}{0.016}                & 0.223                &           &           \\ \cline{1-9}
4                                                                         & \multicolumn{1}{c|}{9375}               & \multicolumn{1}{c|}{12}                  & \multicolumn{1}{c|}{9327}                           & 36                           & \multicolumn{1}{c|}{3.998}                         & \multicolumn{1}{c|}{3.548}                & \multicolumn{1}{c|}{0.033}                & 0.417                &           &           \\ \cline{1-9}
5                                                                         & \multicolumn{1}{c|}{18750}              & \multicolumn{1}{c|}{188}                 & \multicolumn{1}{c|}{18552}                          & 10                           & \multicolumn{1}{c|}{4.989}                         & \multicolumn{1}{c|}{4.288}                & \multicolumn{1}{c|}{0.051}                & 0.650                &           &           \\ \cline{1-9}
6                                                                         & \multicolumn{1}{c|}{15625}              & \multicolumn{1}{c|}{682}                 & \multicolumn{1}{c|}{14942}                          & 1                            & \multicolumn{1}{c|}{5.956}                         & \multicolumn{1}{c|}{4.986}                & \multicolumn{1}{c|}{0.059}                & 0.910                &           &           \\ \cline{1-9}
\textbf{\begin{tabular}[c]{@{}c@{}}R =  13 m\\ d\_min = 4 m\end{tabular}} & \multicolumn{4}{c|}{\textbf{Number of}}                                                                                                                                 & \multicolumn{4}{c|}{\textbf{Average per Scenario}}                                                                                                                &           &           \\ \cline{1-9}
\textbf{\# of UAV}                                                        & \multicolumn{1}{c|}{\textbf{Scenarios}} & \multicolumn{1}{c|}{\textbf{Collisions}} & \multicolumn{1}{c|}{\textbf{Conflict Resolution}} & \textbf{No Conflict Solutions} & \multicolumn{1}{c|}{\textbf{Average Traffic Flow}} & \multicolumn{1}{c|}{\textbf{Path Load 1}} & \multicolumn{1}{c|}{\textbf{Path Load 2}} & \textbf{Path Load 3} &           &           \\ \cline{1-9}
2                                                                         & \multicolumn{1}{c|}{375}                & \multicolumn{1}{c|}{0}                   & \multicolumn{1}{c|}{314}                            & 61                           & \multicolumn{1}{c|}{2.000}                         & \multicolumn{1}{c|}{1.872}                & \multicolumn{1}{c|}{0.032}                & 0.096                &           &           \\ \cline{1-9}
3                                                                         & \multicolumn{1}{c|}{2500}               & \multicolumn{1}{c|}{18}                  & \multicolumn{1}{c|}{2418}                           & 64                           & \multicolumn{1}{c|}{2.993}                         & \multicolumn{1}{c|}{2.640}                & \multicolumn{1}{c|}{0.071}                & 0.282                &           &           \\ \cline{1-9}
4                                                                         & \multicolumn{1}{c|}{9375}               & \multicolumn{1}{c|}{367}                 & \multicolumn{1}{c|}{8972}                           & 36                           & \multicolumn{1}{c|}{3.960}                         & \multicolumn{1}{c|}{3.330}                & \multicolumn{1}{c|}{0.104}                & 0.526                &           &           \\ \cline{1-9}
5                                                                         & \multicolumn{1}{c|}{18750}              & \multicolumn{1}{c|}{2223}                & \multicolumn{1}{c|}{16517}                          & 10                           & \multicolumn{1}{c|}{4.880}                         & \multicolumn{1}{c|}{3.952}                & \multicolumn{1}{c|}{0.124}                & 0.804                &           &           \\ \cline{1-9}
6                                                                         & \multicolumn{1}{c|}{15625}              & \multicolumn{1}{c|}{3984}                & \multicolumn{1}{c|}{11640}                          & 1                            & \multicolumn{1}{c|}{5.738}                         & \multicolumn{1}{c|}{4.517}                & \multicolumn{1}{c|}{0.127}                & 1.095                &           &           \\ \cline{1-9}
\textbf{\begin{tabular}[c]{@{}c@{}}R =  26 m\\ d\_min = 3 m\end{tabular}} & \multicolumn{4}{c|}{\textbf{Number of}}                                                                                                                                 & \multicolumn{4}{c|}{\textbf{Average per Scenario}}                                                                                                                & \textbf{} & \textbf{} \\ \cline{1-9}
\textbf{\# of UAV}                                                        & \multicolumn{1}{c|}{\textbf{Scenarios}} & \multicolumn{1}{c|}{\textbf{Collisions}} & \multicolumn{1}{c|}{\textbf{No Conflict Solutions}} & \textbf{Conflict Resolution} & \multicolumn{1}{c|}{\textbf{Average Traffic Flow}} & \multicolumn{1}{c|}{\textbf{Path Load 1}} & \multicolumn{1}{c|}{\textbf{Path Load 2}} & \textbf{Path Load 3} &           &           \\ \cline{1-9}
2                                                                         & \multicolumn{1}{c|}{375}                & \multicolumn{1}{c|}{0}                   & \multicolumn{1}{c|}{302}                            & 73                           & \multicolumn{1}{c|}{2.000}                         & \multicolumn{1}{c|}{1.946}                & \multicolumn{1}{c|}{0.000}                & 0.053                &           &           \\ \cline{1-9}
3                                                                         & \multicolumn{1}{c|}{2500}               & \multicolumn{1}{c|}{0}                   & \multicolumn{1}{c|}{2430}                           & 70                           & \multicolumn{1}{c|}{3.000}                         & \multicolumn{1}{c|}{2.846}                & \multicolumn{1}{c|}{0.004}                & 0.150                &           &           \\ \cline{1-9}
4                                                                         & \multicolumn{1}{c|}{9375}               & \multicolumn{1}{c|}{2}                   & \multicolumn{1}{c|}{9336}                           & 37                           & \multicolumn{1}{c|}{3.999}                         & \multicolumn{1}{c|}{3.703}                & \multicolumn{1}{c|}{0.015}                & 0.280                &           &           \\ \cline{1-9}
5                                                                         & \multicolumn{1}{c|}{18750}              & \multicolumn{1}{c|}{20}                  & \multicolumn{1}{c|}{18720}                          & 10                           & \multicolumn{1}{c|}{4.998}                         & \multicolumn{1}{c|}{4.525}                & \multicolumn{1}{c|}{0.035}                & 0.437                &           &           \\ \cline{1-9}
6                                                                         & \multicolumn{1}{c|}{15625}              & \multicolumn{1}{c|}{51}                  & \multicolumn{1}{c|}{15573}                          & 1                            & \multicolumn{1}{c|}{5.996}                         & \multicolumn{1}{c|}{5.319}                & \multicolumn{1}{c|}{0.062}                & 0.614                &           &           \\ \cline{1-9}
\textbf{\begin{tabular}[c]{@{}c@{}}R =  26 m\\ d\_min = 5 m\end{tabular}} & \multicolumn{4}{c|}{\textbf{Number of}}                                                                                                                                 & \multicolumn{4}{c|}{\textbf{Average per Episode}}                                                                                                                 & \textbf{} & \textbf{} \\ \cline{1-9}
\textbf{\# of UAV}                                                        & \multicolumn{1}{c|}{\textbf{Scenarios}} & \multicolumn{1}{c|}{\textbf{Collisions}} & \multicolumn{1}{c|}{\textbf{Conflict Resolution}} & \textbf{No Conflict Solutions} & \multicolumn{1}{c|}{\textbf{Average Traffic Flow}} & \multicolumn{1}{c|}{\textbf{Path Load 1}} & \multicolumn{1}{c|}{\textbf{Path Load 2}} & \textbf{Path Load 3} &           &           \\ \cline{1-9}
2                                                                         & \multicolumn{1}{c|}{375}                & \multicolumn{1}{c|}{0}                   & \multicolumn{1}{c|}{314}                            & 61                           & \multicolumn{1}{c|}{2.000}                         & \multicolumn{1}{c|}{1.930}                & \multicolumn{1}{c|}{0.069}                & 0.000                &           &           \\ \cline{1-9}
3                                                                         & \multicolumn{1}{c|}{2500}               & \multicolumn{1}{c|}{0}                   & \multicolumn{1}{c|}{2436}                           & 64                           & \multicolumn{1}{c|}{3.000}                         & \multicolumn{1}{c|}{2.801}                & \multicolumn{1}{c|}{0.194}                & 0.004                &           &           \\ \cline{1-9}
4                                                                         & \multicolumn{1}{c|}{9375}               & \multicolumn{1}{c|}{4}                   & \multicolumn{1}{c|}{9335}                           & 36                           & \multicolumn{1}{c|}{3.999}                         & \multicolumn{1}{c|}{3.622}                & \multicolumn{1}{c|}{0.361}                & 0.016                &           &           \\ \cline{1-9}
5                                                                         & \multicolumn{1}{c|}{18750}              & \multicolumn{1}{c|}{66}                  & \multicolumn{1}{c|}{18674}                          & 10                           & \multicolumn{1}{c|}{4.996}                         & \multicolumn{1}{c|}{4.399}                & \multicolumn{1}{c|}{0.560}                & 0.437                &           &           \\ \cline{1-9}
6                                                                         & \multicolumn{1}{c|}{15625}              & \multicolumn{1}{c|}{230}                 & \multicolumn{1}{c|}{15394}                          & 1                            & \multicolumn{1}{c|}{5.985}                         & \multicolumn{1}{c|}{5.139}                & \multicolumn{1}{c|}{0.787}                & 0.058                &           &           \\ \cline{1-9}
\end{tabular}
\label{tab2}
}
\end{table*}
\normalsize

\begin{figure*}
     \centering
     \subfloat[Travel time for radius-velocity combinations.]{\includegraphics[width=0.33\textwidth]{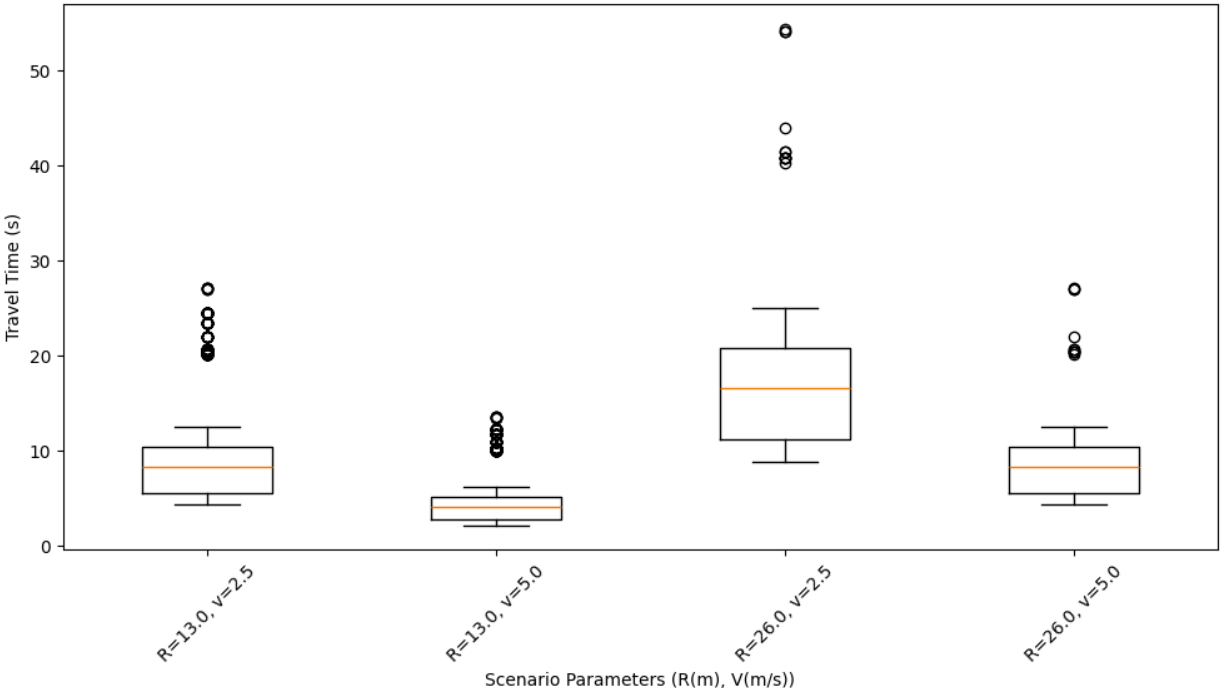}\label{f7}}
     \hfill
     \subfloat[Conflict-pairs with initial displacements.]{\includegraphics[width=0.3\textwidth]{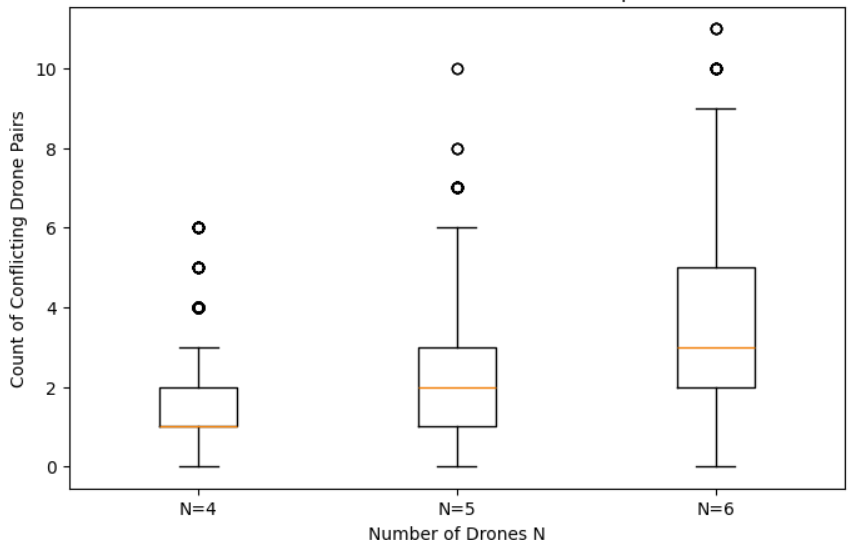}\label{f8}}
     \hfill
     \subfloat[Conflict-pairs with initializations.]{\includegraphics[width=0.3\textwidth]{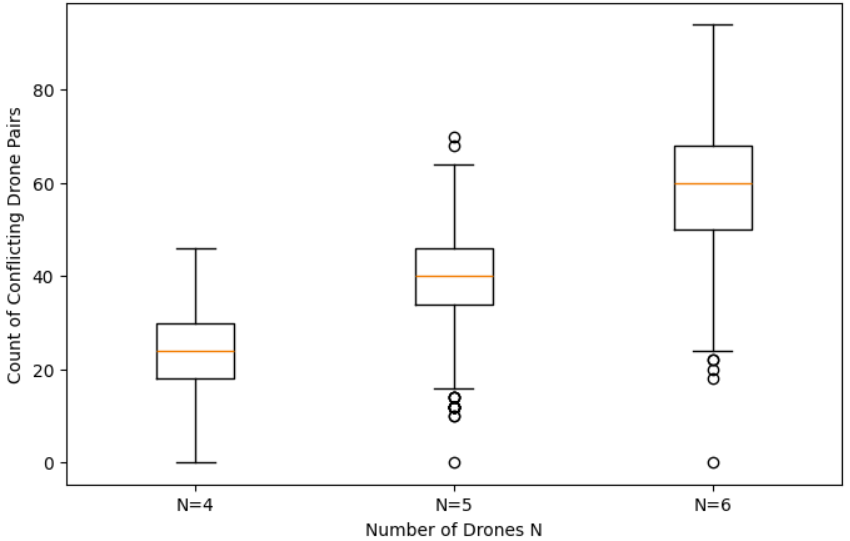}\label{f9}}
     \caption{Output Distributions.}
\end{figure*}

The average traffic flow is optimal for two and three UAVs in most cases, meaning that once they synchronize their entry arrival, the desired traffic flow is carried out. Note that this result is velocity-independent; the only condition is to start simultaneously with any velocity since they are designed geometrically conflict-free for the given minimum clearance. Additionally, the number of collisions increases as this parameter increases, making this design intrinsically infeasible. However, the traffic flow per scenario is near optimal, leading to at least 95\% of paths being conflict-free. Both Sphereabout radius and buffer zone distance influence this parameter. The larger sphere, which can accommodate a wider road, is more effective for collision-free, maximum throughput traffic. On the other hand, the buffer zone distance has an adverse effect, except for the external uncertainties like wind, which is an operational-related design factor, and perhaps can be improved in future. Therefore, we can have conflict-free traffic flow for four UAVs as well, which has at most four unresolved paths. Most of the collisions that happened were due to the overlapping of two left-turn traffic flows. For example, the West-South and North-West flow pair is among the most recurrent collision pairs. 
It can be observed that most UAVs would be assigned to the direct flight path due to optimality. The other paths are less frequent, sharing fewer assignments since arcs are longer. Note that the greater the number of UAV, the more frequent the arc-based path assignment of both kinds, which aligns with the expected optimization solutions. Nonetheless, path three will lose the share in the larger sphere compared to the smaller Sphereabout, while path two almost maintains the path utilization in both spheres. This is potentially due to the longer travel path as the radius increases, leading to a larger distance difference between the optimal one. 

The travel time distribution at the intersection does not seem to depend on the number of UAVs and the minimum safe distance. However, different geometries and velocity distributions can affect the average travel time distribution (Fig. \ref{f7}). It is noticeable that waiting time at the intersection can reach around a minute, which is sufficiently long. If a larger Sphereabout is considered, the velocity must be adjusted to prevent long waiting times. As a result, an optimal design specification can be determined to have conflict-free and lower intersection travel time, ideally a trade-off between radius, minimum separation, and initial velocity.


\subsection{Sensitivity Analysis}

Velocity level can cause a delay in the time between the mobility systems and skip potential conflicts. Therefore, it is evaluated using the following strategies. 

\paragraph{Fixed-lagged Initialization} 

One way is to identify the time difference at the collision and make the positional adjustment as if one UAV is delayed. In other words, a Fixed-lagged distance is added, so the initial position is far ahead or behind such that it takes the other UAV a fixed time difference to reach the collision point if the same velocity is kept. In this regard, first, the time lag for every collision is considered the time difference in their arrival time at the exit point. Afterwards, this time distribution would be shifted backward for one pair collision at the initial time while the velocity is kept constant.     
The corresponding Monte-Carlo simulation with 3000 experiments shows the conflict count distribution for collision set solutions. 
It can be seen from Fig. \ref{f8} that conflict-free paths are achievable for the setting $R=13\;m$ and $d_{min}=3\;m$. The conflict pairs have also become significantly lower, suggesting this approach can be efficient. 

\paragraph{Random Initialization} 

Another way is to adjust the velocity to different random initializations and run a Monte-Carlo simulation, checking for valid and non-collision solutions. In this regard, the velocity is considered as a random number $v \in [1,5] \; \frac{m}{s}$. It is necessary to check whether it is possible to have non-collision paths for all possible traffic flows and UAV entry for this spherical roundabout to ensure this design favours the conflict-aware path planning intersection. To this end, with the same setting, the simulation of the random velocity is shown in Fig. \ref{f9}. 
The collision distribution can be zero for some instances, which suggests the feasibility of collision-free traffic management at Sphereabout for every possible combination and traffic flow. Therefore, velocity modification can be one of the promising alternatives to maximize the traffic throughput.

\section{Conclusion}\label{conc}

The work addresses the practical and operational three-dimensions of traffic management for networked urban air mobility systems. First, the \emph{Sphereabout} concept is introduced as the novel spherical intersection design. This setup incorporates the symmetrical geometry of the intersection and the safety metrics, leading to integration of tube corridors with the central sphere for conflict-free traffic flows. The spherical structure can potentially resolve collisions by time-independent path assignment and a time-dependent technique, such as lagging displacement, for all feasible traffic flows. This design can help the urban aerial traffic control units, since it is simple and inspired by the urban environment network and guarantees intrinsic geodesic-inspired traffic management for less than three operational UAV and, with further improvement, can cover up to four.  

There are several directions in this research for further consideration. First, the optimization framework should incorporate bi-objective non-linear model of conflict resolution while minimizing the waiting time and energy consumption at the intersection. More advanced and dynamic path planning is worth exploring in aerial network, especially at the scale of spherical roundabout, which requires more structured and risk-free transition between entries, however with sufficient velocity to avoid hovering and time delay. Finally, the higher congestion scenarios handling dynamic demand at the network level are another avenue which will be considered in future directions.  




\end{document}